\documentclass[11pt]{article}
\usepackage{amssymb}
\usepackage{amsmath}
\usepackage{amsfonts}
\usepackage{amsthm}

\numberwithin{equation}{section}

\newtheorem{theorem}{Theorem}[section]
\newtheorem{lemma}{Lemma}[section]
\newtheorem{proposition}{Proposition}[section]

\newtheorem{remark}{Remark}[section]
\newcommand\re[1]{(\ref{#1})}
\def\R{\mathbb{R}}

\def\T{\mathbb{T}}
\def\eps{\varepsilon}

\begin{document}
\title{A note on ill-posedness   for the KdV equation}
\author{\normalsize \bf  Luc Molinet  \\
{\footnotesize \it L.M.P.T., Universit\'e Fran\c cois Rabelais Tours, F\'ed\'eration Denis Poisson-CNRS,} \\
{\footnotesize \it Parc Grandmont, 37200 Tours,  France} }
\maketitle
\begin{abstract}
We prove that the solution-map $ u_0 \mapsto u $ associated with the KdV equation cannot be continuously extended
 in $ H^s(\R) $ for $ s<-1 $.  The main ingredients are the well-known Kato smoothing effect for the mKdV equation as well as  the Miura transform.
  \end{abstract}
\section{Statement of the result}
  The KdV equation is a canonical dispersive equation that reads
   \begin{equation}\label{KdV}
  v_t+v_{xxx} -6 v v_x =0 \; 
  \end{equation}
  where $ u(t,x) $ is a real valued function. 
   In the end of the sixties, Miura discovered that the KdV equation is related to the defocussing mKdV equation
    \begin{equation}\label{mKdV}
    u_t+u_{xxx} -  6 u^2 u_x =0
  \end{equation}
   by the now so-called   Miura transformation. Indeed, the Miura transformation $ \Phi : u\mapsto u_{x}+u^2 $  maps  a smooth solution of the defocussing mKdV  to a solution of the real-valued KdV equation

    The Cauchy problem,  in Sobolev spaces,  associated with these equations has been extensively studied since the end of the eighties (see for instance \cite{KPVCPAM93}, \cite{Bourgain1993}, \cite{KPVJAMS96} and references herein). The lowest Sobolev index reached for the well-posedness theory is 
     $ s_u(KdV)=-3/4 $ for the KdV equation (cf. \cite{Kis} or \cite{Guo2009}) and $s_u(mKdV)=1/4 $ for the mKdV equations (cf. \cite{KPVCPAM93}). Note that the difference of one derivative between these two results agrees with the  lost of one derivative in the Miura transformation. These indexes have been proved to be optimal if one requires moreover the solution-map  $ u_0\mapsto u $ to be uniformly continuous on
     bounded set  from $ H^s(\R) $ into $ C([0,T];H^s(\R)) $ (cf. \cite{KPVDuke01}, \cite{CCT2003}). 
     However, one may expect that the solution-map associated with these equations can be
      continuously extended below these indexes especially since  they are far above  the critical Sobolev indexes for the dilation symmetry that  are
       respectively $ s_{d}(KdV)=-3/2 $ and $s_{d}(mKdV)=-1/2 $.
      It is worth emphasizing that such results have been proved in the periodic case by Kappeler and Topalov (cf. 
       \cite{KT1}, \cite{KT2}). More precisely, they proved that the solution-maps associated to the KdV and mKdV equations  can be uniquely continuously extended in respectively $ H^{-1}(\T) $ and $ L^2(\T) $.
       Note that, in the periodic setting, these  solution-maps are known to be non uniformly continuous on
      bounded sets of $ H^s(\T) $ for respectively $ s<-1/2 $ and $ s<1/2 $. \\

      In this note we combine the Miura transforms and the so-called local Kato smoothing effect to prove that the solution-map associated with the KdV equation cannot be continuously extended in $ H^s(\R ) $ for $ s<-1 $.
           \begin{theorem}\label{maintheorem}
         Let  $ s<-1$ be given. For any $ \varepsilon>0 $  there exists $ h_{\varepsilon} \in H^\infty(\R) $ with 
          $ \| h_{\varepsilon}\|_{L^2} =\varepsilon $ such that, for any $ T>0 $, the solution-map $ u_0 \mapsto u $ associated with the  KdV equation is discontinuous at $ h_\varepsilon $ from  $ H^\infty(\R) $  equipped with
   the $ H^s(\R) $ topology into  $ {\mathcal D}'(]0,T[\times \R) $.
 \end{theorem}
 \begin{remark}\label{remark}
 Let  $ B_{H^{-1}}(0,R) $ be the ball of $ H^{-1}(\R) $ centered at the origin with radius $ R$. Actually we also prove that there exits no $ R>0 $ and no $ T>0 $ such that the flow-map $ u_0\mapsto u $ associated with the KdV equation  is continuous from 
  $ B_{H^{-1}}(0,R) \cap H^\infty(\R) $, equipped with the weak topology of $ H^{-1}(\R) $, into $ {\mathcal D}'(]0,T[\times \R) $.
 \end{remark}
\section{Proof.} 
The idea of the proof is the following:   On one hand the Miura transformation is  discontinuous in $ H^s(\R) $ as soon as $ s<0 $ so that  we can find a sequence of initial data , bounded in $ L^2(\R) $, that converges to $0 $ in $ H^s(\R) $ but 
 which  Miura transform converges to some non identically vanishing element $ \theta $ of $ H^\infty(\R)$.  
On the other hand, using the Kato smoothing effect, we can prove that the sequence of associated  solutions  $ \{u_N \} $ to the mKdV equation converges, up to a subsequence, to some   solution 
   $ u $  of the mKdV equation that belongs to $C_{w}(\R;H^{-1}(\R)) $ and such that $ u(0)=0 $. This forces the sequence of associated Miura transform $ \{\Phi(u_{N_k}) \} $ to converge towards $ \Phi(u) $ in some sense.  This yields the result since
   we  can prove that $ \{\Phi(u) \} $ is not the classical solution of the KdV equation emanating from $ \theta$. \vspace*{2mm}

\noindent
{\it Step 1- Choice of the sequence of initial data -}   Let $ \{h_{N}\} \subset H^\infty(\R) $ be the sequence defined by
$$
\hat{h}_{N}(\xi) = \chi_{ [N,N+1] }(|\xi|) \, .
$$
Clearly $ \{h_{N}\} $ is bounded in $L^2(\R) $. It tends to $ 0 $ in $ H^s(\R) $ for any  $ s<0 $ and weakly in $ L^2(\R) $.
Let us now fix $ 0<\varepsilon<1$ and consider the sequence of initial data $ \{\varepsilon h_N\} $. Straightforward calculations
      lead to
     $$
     \Phi(\eps h_N)= \eps^2 \theta + \eps  h_N' +\eps^2 \alpha_N
     $$
     where
     $$
     \hat{\theta}(\xi)=(1-|\xi|) \, \chi_{[0,1]}(|\xi|) \quad \Bigl(i.e.   \; \theta(x)=\frac{2}{x^2} (1-\cos (x)) , \quad \forall x\in \R \Bigr)
     $$
     and
     $$
     \hat{\alpha}_N(\xi) =\Bigl(1-\Bigl||\xi|-2N-1\Bigr|\Bigr)\,\chi_{[0,1]}( ||\xi|-2N-1|) \, .
     $$
     Clearly  $ \theta \in H^\infty(\R) $ with $ \|\theta\|_{L^2}\sim 1$, 
     $ \alpha_N \to 0 $ in $ H^{-1}(\R)$ and
       $  h_N' \to 0 $ in $ H^{s}(\R)$ for $s<-1 $.  In particular, for any $ s<-1$,
   \begin{equation} \label{gu}
   \Phi(\eps h_N)\to \eps^2\theta \mbox{ in } H^s(\R)  \mbox{ as } N\to +\infty\, .
   \end{equation}
    {\it Step 2 - Uniform bounds on the sequence of  emanating solutions -}
   Since $ h_{N}\in H^\infty(\R) $, it follows from classical well-posedness results  for mKdV (see for instance \cite{GTV}) that the solution $ u_N$ of \eqref{mKdV}  emanating from $ \varepsilon h_N $  exists for all times and belongs to
 $ C(\R;H^\infty(\R)) $. Moreover, it is well known that the $ L^2$-norm of the solution is a constant of
  the motion and that  the so-called Kato smoothing effect holds. For sake of completeness
  we give hereafter a version of  this smoothing effect that is suitable for our purpose
  (see for instance
  \cite{GTV} for a general setting) and use it to prove some uniform continuity result.\vspace{2mm} \\
 \begin{lemma} \label{PropositionKato}Let $ u_0\in  H^\infty(\R) $. Then for any $ T>0$ and $ R>0 $ there exists $ C(T,R) >0 $ such that the  emanating solution of
  \re{mKdV} safisfies
  \begin{equation}\label{p1}
  \int_{-T}^T \int^{-R}_R |\partial_x u(t,x)|^2 dx\, dt \le C(R,T) \|u_0\|_{L^2}^2 \; .
  \end{equation}
  Moreover,  for any function $ \varphi\in C^\infty_c(\R) $ with compact support in $ ]-R,R[ $ and any interval 
   $ [t,t+\delta] \subset ]-T,T[ $ with $ 0<\delta<1$,  there exists $ C(\varphi, T,R) >0 $ such that 
   \begin{equation}\label{p2}
   \Bigl|   \int_{\R}\varphi u(t+\delta )\, dx - \int_{\R} \varphi u(t)\, dx \Bigr|\le C(\varphi,R,T)\, \delta^{1/2} \Bigl(\|u_0\|_{L^2}^2
    + \|u_0\|_{L^2}^3\Bigr)\; .
  \end{equation}
   \end{lemma}
  \begin{proof}
  Let $ \phi\in C^\infty_c(\R) $  with $ \phi\equiv 1 $ on $[-1,1]$, $ 0\le \phi \le 1 $ on $ [-2,2] $ and 
   $\mbox{Supp } \phi\subset ]-2,2[$.  For $ R>1 $
   we set $ h'(x):=\phi(x/R) $ and $ h(x):=\int_{-\infty}^x h'(t)\, dt $. Integrating \re{mKdV} against $ h u $ one obtains
   $$
   \frac{d}{dt} \int_{\R}h u^2 +3 \int_{\R} h' (u_x)^2 -\int_{\R} h^{'''} u^2 +3  \int_{\R} h' u^4 = 0
   $$
    Since the $ L^2 $-norm is conserved by the flow of \re{mKdV}, we get
    $$
   \frac{d}{dt} \int_{\R}h u^2 + \int_{\R} h' (u_x)^2 \le C \,  \|u_0\|^2_{L^2} \; .
    $$
    Integrating this inequality in $(-T,T) $ and using the properties of $ h'$, we thus infer that
    $$
     \Bigl| \int_{-T}^T  \int_{-R}^R  (u_x)^2 \Bigr|\le  \Bigl| \int_{-T}^T \int_{\R}  h' (u_x)^2 \Bigr|\le C\,  (T+R) \|u_0\|^2_{L^2}\; .
    $$
    This completes the proof of \eqref{p1}. To prove \eqref{p2} we integrate \re{mKdV} against $ \varphi \in C^\infty_c(\R) $,
     compactly supported in $ ]-R,R[ $,   to get 
  $$
    \frac{d}{dt} \int_{\R}\varphi u -\int_{\R} \varphi_{xxx} u  + 2\int_{\R} \varphi' u^3 = 0\; .
  $$
    Integrating this identity on $ [t,t+\delta ] \subset (-T,T) $ and using H\"older inequality in space, we obtain
   \begin{equation}\label{ez}
  \Bigl|   \int_{\R}\varphi u(t+\delta )\, dx - \int_{\R} \varphi u(t)\, dx \Bigr|\le C(\varphi,T)\int_{t}^{t+\delta}  \Bigl[\Bigl(\int_{-R}^{R}  u^2\, dx\Bigr)^{1/2}+
  \Bigl(\int_{-R}^R |u|^6\,dx\Bigr)^{1/2} \Bigr]
     \end{equation}
   On the other hand, interpolating between  \eqref{p1} and the conservation of the $ L^2 $-norm, we infer by Sobolev inequality that  
     \begin{equation}\label{ezz}
   \int_{-T}^T \int_{-R}^R |u|^6 \, dx\, dt \le C(T,R) \| u_0\|_{L^2}^6Ê\; .
     \end{equation}
   \eqref{p2} then follows from  \eqref{ez}  and \eqref{ezz} by applying H\"older inequality in time.
  \end{proof}
\noindent
{\it Step 3 - Convergence results and properties of the limit -} From the above lemma we deduce the following convergence result :
\begin{proposition}
There exists a subsequence   $ \{u_{N_k}\} $ of $ \{u_{N}\} $ such that
 \begin{align}
        u_{N_k} \to u \quad & \mbox{a.e. in } \R^2,\label{pro1} \\
         u_{N_k} \rightharpoonup  u  \quad & \mbox{ weakly star in }  L^\infty(\R;L^2(\R)) ,\label{pro2}\\
         u_{N_k}^2 \rightharpoonup u^2 \quad & \mbox{ weakly star  in } L^\infty(\R;H^{-1}(\R)),\label{pro3}
        \end{align}
where $u \in C_{w}(\R;L^2(\R)) \cap L^2_{\rm loc}(\R; H^1_{\rm loc}) $ satisfies the mKdV equation in the distributional sense with $ u(0)=0 $.
\end{proposition}
  \begin{proof}
We first notice that,  by the conservation of the $ L^2 $-norm, 
       $ \{ u_N\} $ is bounded in $ L^\infty(\R;L^2(\R)) $. Second, from  \eqref{p1} in Lemma  \ref{PropositionKato}, $ \{\partial_x u_N\} $ is bounded in $ L^2_{\rm loc}( \R^2) $.  Third, by the equation,
  $ \{\partial_{t}u_{N}\} $  is bounded in
        $ L^2_{\rm loc}(\R;H^{-2}_{\rm loc}) $.
 It thus follows from Aubin-Lions compactness theorem that there exist
         $ u\in L^\infty(\R; L^2(\R))\cap L^2_{\rm loc}(\R; H^{1}_{\rm loc}) $  and an  increasing sequence of integers $ \{N_k\}$ such that \eqref{pro1}-\eqref{pro3} hold. Now, using \eqref{pro1}, we can pass to the limit on the equation to infer that $ u $ satisfies the mKdV equation in the distributional sense. Moreover, it follows from \eqref{p2} that, for any fixed
 $ T>0 $ and any fixed         $\varphi\in C^\infty_c(\R) $, the sequence of functions $\{ t\mapsto \int_{\R} \varphi u_{N_k} (t) \, dx \}_{k\ge 0}$ is uniformly equi-continuous on $ [-T,T] $ and thus converges in $ C(-T,T) $ towards some continuous function $ w $  thanks to  Ascoli's theorem. \eqref{pro1} then ensures that $ w= \int_{\R} \varphi u(t) \, dx  $ which proves that  $ u \in C_{w}(\R;L^2(\R)) $. Finally,  since $
          u_{N}(0)=\varepsilon h_N \rightharpoonup 0 $ in $L^2(\R) $ we infer that $ u(0)=0 $.
\end{proof}
\begin{proposition}\label{prop2}
The sequence of solutions  $ \{\Phi(u_{N_k})\} $ of the KdV equation converges  weakly star in $ L^\infty(\R;H^{-1}(\R))$ 
 towards $ \Phi(u) $  that satisfies : For all $ \varphi \in C^\infty_c(\R) $, $t\mapsto  \int_{\R} \varphi(x) \Phi(u(t,x))\, dx $ is, 
   up to modifications on a set of measure zero, a continuous function on $ \R $ with value zero at the origin.
   \end{proposition}
\begin{proof}
First from \eqref{pro2}-\eqref{pro3} we infer that, for any $ T>0$, 
\begin{equation}
\Phi(u_{N_k}) \rightharpoonup \Phi(u) \mbox{ weakly star  in } L^\infty(\R;H^{-1}(\R)) \label{pro4}
\end{equation}
and thus in $ D'(\R^2) $. Now, to prove the remaining of the statement, we use an exterior regularization by a  mollifier
sequence $ \{\rho_m\} $ where 
$$
\rho_m(x) =\Bigl( \int_{\R} \rho(y) \, dy \Bigr)^{-1} m \rho(mx), \quad x\in\R, \; n\ge 1 ,
$$
and  $ \rho \not \equiv 0 $ is a non negative  $ C^\infty $ function compactly supported in $ ]-1,1[ $. We set $ u_m=\rho_m\star u $. 
Since  $ u\in L^\infty(\R; L^2(\R))\cap L^2_{\rm loc}(\R; H^{1}_{\rm loc}) $ and  satisfies the mKdV equation in the distributional sense, one can easily check that    $ u_t\in L^2_{\rm loc}(\R; H^{-2}_{\rm loc})$. It follows that 
$$
\partial_t u_m +\partial_{x}^3 u_{m} = 2  \rho_m' \ast (u^3) \; 
$$
 with $ \partial_t u_m\in L^2_{\rm loc} (\R^2) $. Multiplying this equation with $\varphi u_m $ and integrating by parts on $ \R $ we get 
\begin{align*}
\frac{d}{dt}\Bigl(  \int_{\R} \varphi u_m^2 \, dx\Bigr) +& 3\int_{\R} \varphi' (\partial_x u_{m})^2  -\int_{\R} \varphi^{'''} u_m^2
+3Ê\int_{\R} \varphi' u_m^4 \\
&=-4 \int_{\R} (\varphi' u_m+ \varphi \partial_x u) [u_m^3- \rho_m\ast u^3] \end{align*}
Integrating in time between $ 0 $ and $ t>0 $, recalling that $ u(0,\cdot)=0$,  this leads to 
\begin{align}
 \int_{\R} \varphi u_m^2(t) \, dx +& 3\int_0^t \int_{\R} \varphi' (\partial_x u_{m})^2  -\int_0^t\int_{\R} \varphi^{'''} u_m^2
+3\int_0^tÊ\int_{\R} \varphi' u_m^4 \nonumber \\
&=-4 \int_0^t\int_{\R} (\varphi' u_m+ \varphi \partial_x u_n) [u_m^3- \rho_m\ast u^3] \, .\label{mm}\end{align}
Now, since  $ u\in L^\infty(\R; L^2(\R))\cap L^2_{\rm loc}(\R; H^{1}_{\rm loc}) $ we infer that 
$$
\partial_x u_m \to u_x \mbox{ in } L^2_{\rm loc}(\R^2), \; u_m\to u \mbox{ in }   L^6_{\rm loc}(\R^2)\;  \mbox{Êand }
Ê\rho_m\ast u^3 \to u^3  \mbox{ in } L^2_{\rm loc}(\R^2)\, .$$
Moreover, there exists a mesurable set $ E\subset \R $ with $ \mbox{mes }(E)=0 $ such that 
 $ u(t) \in L^2(\R) $ for all $ t\in \R/E $.  For any $ t\in \R/E $ this ensures that  $ u_m(t) \to u(t) $ in $L^2(\R) $ and thus passing to the limit in \eqref{mm}  as $ m $ goes to infinity we obtain
\begin{equation}
 \int_{\R} \varphi u^2(t) \, dx=  -3\int_0^t \int_{\R} \varphi' (\partial_x u)^2  +\int_0^t\int_{\R} \varphi^{'''} u^2
-3\int_0^tÊ\int_{\R} \varphi' u^4 \, .\label{ppp}
\end{equation}
This shows that, up to a modification on $ E $, $ t\mapsto   \int_{\R} \varphi u^2(t) \, dx$ is a continuous function with value zero at the origin.

Recalling that  $ u\in C_{w}(\R;L^2(\R)) $ with $ u(0)=0$ and that $ \Phi(u)=u'+u^2 $, this completes the proof of the proposition.
\end{proof}

\noindent
 {\it - Conclusion -} Let us fixed $s<-1$ and $ \varepsilon>0$. From  \re{gu},  
 $\Phi(\varepsilon h_N)\to\varepsilon^2 \theta\neq 0 $ in $ H^{s}(\R) $, with $\|\theta\|_{L^2} \sim 1 $,  and  the solution $v $
   to \eqref{KdV} emanating from $  \eps^2 \theta$  belongs to $C(\R;H^{\infty}(\R)) $. Let $ \phi\in C^\infty_c(\R) $  with $ \phi\equiv 1 $ on $[-1,1]$, $ 0\le \phi \le 1 $ on $ [-2,2] $ and 
   $\mbox{Supp } \phi\subset ]-2,2[$. Setting $ \varphi:=\phi \theta\in C^\infty_c(\R)  $
    we thus infer that there exists $ r_0>0 $ such that 
   \begin{equation}\label{conc1}
   \int_{\R} \varphi(x) v(t,x) \, dx \ge\frac{\varepsilon^2}{2}  \int_{-1}^{1}|\theta(x)|^2 \, dx>0  \mbox{ for any } t\in ]-r_0,r_0[ \; .
   \end{equation}
   On the other hand, according to  Proposition \ref{prop2}, there exists $ u\in L^\infty(\R; L^2(\R))\cap L^2_{\rm loc}(\R; H^{1}_{\rm loc}) $ and an  increasing sequence of integers $ \{N_k\}$ such that $ u_{N_k}  \rightharpoonup u $ weakly star  in $ L^\infty(\R;L^2(\R))$ and  such that  the sequence of solutions $\{ v_{N_k}\}:= \{\Phi(u_{N_k})\}$ of KdV, emanating from
  $ \{ \Phi(\varepsilon h_N)\} $, tends weakly  star towards $ \Phi(u)$ in  $L^\infty(\R; H^{-1}(\R)) $. Moreover, there exists 
    $ 0<r<r_0 $ such that 
       \begin{equation}\label{conc2}
       \Bigr| \int_{\R} \varphi(x) \Phi(u(t,x))\, dx \Bigl|\le \frac{\varepsilon^2}{4}  \int_{-1}^{1}|\theta(x)|^2 \, dx \mbox{ for a.e. } t\in ]0,r[ \; . 
       \end{equation}
Gathering \eqref{conc1} and \eqref{conc2} this ensures that for almost every $ t\in ]0,r[$, 
       $$
        v(t) \neq \Phi(u(t)) \mbox{ a.e. in } \R \; .
        $$
        which completes the proof of the theorem. Finally note that Remark \ref{remark} follows from the fact that 
         $ \|\Phi(\varepsilon h_N) \|_{H^{-1}(\R)} \sim \varepsilon $, for all $ n\ge 1$,   and $ \Phi(\varepsilon h_N) \rightharpoonup \varepsilon^2 \theta $ in $ H^{-1}(\R) $.      \vspace{2mm} \\
   \noindent
 {\bf Acknowledgements:}   L.M.  was partially supported by the ANR project
 "Equa-Disp". He  would like to thank Didier Pilod for his careful reading of a preliminary version of this work.

  \end{document}